\documentstyle[12pt]{article}
\newcommand{\nt}{\noindent}
\newcommand{\oid}[3]{\mbox{${\cal #1}_{#2}^{#3}$}}
\newcommand{\idn}[3]{\mbox{${\bf #1}_{#2}^{#3}$}}
\newcommand{\idrel}[3]{\mbox{${\cal #1} \stackrel{1}{#2} {\cal #3}$}}

\begin{document}
\newtheorem{theorem}{Theorem}[section]
\newtheorem{prop}{Proposition}[section]
\newtheorem{lemma}{Lemma}[section]
\newtheorem{cor}{Corollary}[section]
\newtheorem{rem}{Remark}[section]

\title{Operators with extension property and the principle of local 
reflexivity}
\author
{Frank Oertel\\
Zurich}
\maketitle

\nt{\it{Key words and phrases:}} geometry of Banach spaces, operator ideals, 
tensor norms, accessibility, principle of local reflexivity\\ 

\begin{abstract}
{\nt}Given an arbitrary $p$-Banach ideal (\oid{A}{}{}, \idn{A}{}{}) ($0 < p \leq 1$), 
we ask for geometrical properties  of (\oid{A}{}{}, \idn{A}{}{}) which are sufficient 
(and necessary) to allow a transfer of the principle of local reflexivity to 
(\oid{A}{}{}, \idn{A}{}{}).\\

{\nt}{\it{1991 AMS Mathematics Subject Classification:}} primary 46M05, 47D50; 
secondary 47A80
\end{abstract}

\section{Introduction}

{\nt}Given Banach spaces $E$, $F$ and a maximal Banach ideal (\oid{A}{}{}, \idn{A}{}{}),
we are interested in reasonable sufficient conditions on $E$, $F$ and
(\oid{A}{}{}, \idn{A}{}{}) such that (\oid{A}{}{}, \idn{A}{}{})
is accessible. In general it is a nontrivial subject to prove accessibility of  
maximal Banach ideals since non-accessibility only appears on Banach
spaces without the metric approximation property, and in 1991, {\sc{Pisier}}
made use of such a Banach space (the Pisier space $P$) to construct a non-accessible
maximal Banach ideal (cf. \cite{df}, 31.6). On the other hand, accessible Banach ideals
allow a suggestive calculus which leads to further results concerning
the local structure of operator ideals (e.g., a transfer of the principle of
local reflexivity from the operator norm to suitable ideal norms 
\idn{A}{}{} (cf. \cite{df}, \cite{lr}, \cite{oe1}, \cite{oe2} and \cite{oe4})).\\

{\nt}This paper is mainly devoted to the investigation of a class of 
(maximal) Banach ideals which do allow the transfer of the operator 
norm estimation in the principle of local reflexivity to the norm
of the given operator ideal.\\ 

{\nt}We only deal with Banach spaces and most of our notations and definitions 
concerning Banach spaces and operator ideals are standard and can be
found in the detailed monographs \cite{df} and \cite{p1}. However, if 
$(\oid{A}{}{}, \idn{A}{}{})$ and $(\oid{B}{}{}, \idn{B}{}{})$ are given
quasi-Banach ideals,
we will use the shorter notation $(\oid{A}{}{d}, \idn{A}{}{d})$ for the dual
ideal (instead of $(\oid{A}{}{dual}, \idn{A}{}{dual})$) and the abbreviation 
\idrel{A}{=}{B} for the isometric identity
$(\oid{A}{}{}, \idn{A}{}{}) = (\oid{B}{}{}, \idn{B}{}{})$. The inclusion 
$(\oid{A}{}{}, \idn{A}{}{}) \subseteq (\oid{B}{}{}, \idn{B}{}{})$ is often
shortened by \idrel{A}{\subseteq}{B}, and if $T : E \longrightarrow F$
is an operator, we indicate that it is a metric injection by writing 
$T : E \stackrel{1}{\hookrightarrow} F$. Each section of this paper 
includes the more special terminology which is not so common to specialists
in the geometry of Banach spaces. 

\section{Minimal and maximal Banach ideals and their conjugates}
 
{\noindent}Let $E,F$ be arbitrary Banach spaces and $T \in \oid{L}{}{}(E,F)$.
Given a $p$-Banach ideal (\oid{A}{}{}, \idn{A}{}{}) ($0 < p \leq 1$) and a 
$q$-Banach ideal (\oid{B}{}{}, \idn{B}{}{}) ($0 < q \leq 1$), we can construct
further ideals:
\begin{itemize}
\item    $T$ belongs to  $\oid{A}{}{} \circ \oid{B}{}{-1}(E,F)$ if $TS \in 
         \oid{A}{}{}(G, F)$ for all operators $S \in \oid{B}{}{}(G, E)$ and
         for all Banach spaces $G$. Let $\idn{A}{}{} \circ \idn{B}{}{-1}(T) 
         : = \sup \{\idn{A}{}{}(TS): S \in \oid{B}{}{}(G, E), \idn{B}{}{}(S) 
         = 1\}$. Then $(\oid{A}{}{} \circ \oid{B}{}{-1}, \idn{A}{}{} \circ \idn{B}{}{-1})$
         defines a $p$-Banach ideal, called the {\it{right-\oid{B}{}{}-quotient of
         \oid{A}{}{}}}.
\item    $T$ belongs to  $ \oid{B}{}{-1} \circ \oid{A}{}{}(E,F)$ if $ST \in 
         \oid{A}{}{}(E, G)$ for all operators $S \in \oid{B}{}{}(F, G)$ and
         for all Banach spaces $G$. Let $\idn{B}{}{-1} \circ \idn{A}{}{}(T) 
         : = \sup \{\idn{A}{}{}(ST): S \in \oid{B}{}{}(F, G), \idn{B}{}{}(S) 
         = 1\}$. Then $(\oid{B}{}{-1} \circ \oid{A}{}{}, \idn{B}{}{-1} \circ 
         \idn{A}{}{})$ defines a $p$-Banach ideal, called the {\it{left-\oid{B}{}{}-quotient of
         \oid{A}{}{}}}.
\item    $T$ belongs to  $\oid{A}{}{} \circ \oid{B}{}{}(E,F)$ if   
         there exists a Banach space $G$ and operators $R \in \oid{B}{}{}(E,G)$ 
         and $S \in \oid{A}{}{}(G, F)$ such that $T = SR$. If we set 
         $\idn{A \circ B}{}{}(T) : = \inf \{ \idn{A}{}{}(S) \cdot 
         \idn{B}{}{}(R)\}$, the infimum being taken over all possible 
         factorizations of $T$, then 
         $(\oid{A}{}{} \circ \oid{B}{}{}, \idn{A}{}{} \circ \idn{B}{}{})$ is 
         a $r$-Banach ideal, where $\frac{1}{r} =  \frac{1}{p}+\frac{1}{q}$.
\end{itemize} 

{\noindent}Note, that in general  $(\oid{B}{}{-1} \circ \oid{A}{}{}, \idn{B}{}{-1} 
\circ \idn{A}{}{}) \not= (\oid{A}{}{} \circ 
\oid{B}{}{-1}, \idn{A}{}{} \circ \idn{B}{}{-1})$. Using  
{\sc{Pisier's}} counterexample of a non-accessible maximal Banach ideal, it follows  
that there even exists a Banach ideal (\oid{A}{}{}, \idn{A}{}{}) such that  
$(\oid{A}{}{-1} \circ \oid{I}{}{}, \idn{A}{}{-1} \circ \idn{I}{}{}) 
\not= (\oid{I}{}{} \circ \oid{A}{}{-1}, \idn{I}{}{} \circ 
\idn{A}{}{-1})$ (cf. \cite{oe4}), where (\oid{I}{}{},\idn{I}{}{}) denotes the
class of all integral operators.\\ 
 
{\nt}Using the class of all approximable operators 
$(\bar{\oid{F}{}{}} , \| \cdot \|)$, important special cases of the 
previous constructions are given by: 
\begin{itemize}
\item[-]  the {\it{minimal kernel}} of (\oid{A}{}{}, \idn{A}{}{}):
\begin{center}
       $(\oid{A}{}{min}, \idn{A}{}{min}) : = (\bar{\oid{F}{}{}} 
       \circ \oid{A}{}{}\circ \bar{\oid{F}{}{}},\/ 
       \| \cdot\| \circ \idn{A}{}{} \circ \| \cdot\|)$
\end{center}
\end{itemize}

\begin{itemize}       
\item[-]  the {\it{maximal hull}} of (\oid{A}{}{}, \idn{A}{}{}):
\begin{center}
       $(\oid{A}{}{max}, \idn{A}{}{max}) : = (\bar{\oid{F}{}{}}^{-1} 
       \circ \oid{A}{}{}\circ \bar{\oid{F}{}{}}^{-1}, \/
       {\| \cdot\|}^{-1} \circ \idn{A}{}{} \circ {\| \cdot\|}^{-1})$
\end{center}       
\end{itemize}

{\noindent}(\oid{A}{}{}, \idn{A}{}{}) is said to be {\it{minimal}} if 
(\oid{A}{}{}, \idn{A}{}{}) = (\oid{A}{}{min}, \idn{A}{}{min}). If 
(\oid{A}{}{}, \idn{A}{}{}) = (\oid{A}{}{max}, \idn{A}{}{max}), we say that 
(\oid{A}{}{}, \idn{A}{}{}) is {\it{maximal}}. Obviously (\oid{A}{}{min}, 
\idn{A}{}{min}) is the largest minimal operator ideal which is contained in 
(\oid{A}{}{}, \idn{A}{}{}), and (\oid{A}{}{max}, \idn{A}{}{max}) is the 
smallest maximal operator ideal which contains (\oid{A}{}{}, \idn{A}{}{}). 
Although the product of Banach operator ideals need not to be 
normed again, it can be shown that if (\oid{A}{}{}, \idn{A}{}{}) 
is a Banach ideal then (\oid{A}{}{min}, \idn{A}{}{min}) is also a 
{\it{Banach}} ideal (cf. \cite{br}).\\

{\noindent}Concerning a deeper investigation of local properties of operator 
ideals, two further important {\it{Banach}} ideals play a
key role. As before, let (\oid{A}{}{}, \idn{A}{}{}) be a given 
$p$-Banach ideal $(0 < p \leq 1)$, $E, F$ be arbitrary Banach spaces and 
$T \in \oid{L}{}{}(E,F)$ (The notion of the conjugate of an operator 
ideal was introduced by {\sc{Gordon, Lewis, Retherford}}  
(cf. \cite{glr}, \cite{jo})).

\begin{itemize}
\item $T \in \oid{A}{}{\displaystyle\ast}(E,F),$
      if there exists a constant $c \geq 0$ such that for all Banach spaces 
      $E_0, F_0$ and for all operators $B \in \oid{F}{}{}(E_0, E), S\in 
      \oid{A}{}{}(F_0, E_0), A \in \oid{F}{}{}(F, F_0)$ 
      \begin{center}
      $\mid tr(TBSA) \mid \leq c \cdot \|B\| \cdot \idn{A}{}{}(S) \cdot \|A\|$.
      \end{center} 
      {\noindent}Setting
      \begin{center}
      $\idn{A}{}{\displaystyle\ast}(T) : = \inf(c)$ \/, 
      \end{center}
      where the infimum is taken over all such constants $c$, we obtain a Banach 
      ideal \\(\oid{A}{}{\displaystyle\ast}, \idn{A}{}{\displaystyle\ast}), the 
      {\it{adjoint}} of (\oid{A}{}{}, \idn{A}{}{}).

\item  $T \in \oid{A}{}{\Delta}(E, F),$ if there exists a constant $c \geq 0$
       such that for all {\it{finite}} operators \\$L \in {\cal F}(F, E)$
       \begin{center}
       $\mid tr( TL) \mid \leq c \cdot \idn{A}{}{}(L).$
       \end{center}
       Setting 
       \begin{center}
       $\idn{A}{}{\Delta}(T) : = \inf(c)$ \/,
       \end{center}
       where the infimum is taken over all such constants $c$, we obtain 
       a Banach ideal \\(\oid{A}{}{\Delta}, \idn{A}{}{\Delta}), the 
       {\it{conjugate}} of (\oid{A}{}{}, \idn{A}{}{}). 
\end{itemize}

\section{On tensor norms and associated Banach ideals}
First we recall the basic notions of Grothendieck's metric theory of tensor products
(cf., e.g., \cite{df}, \cite{gl}, \cite{gr}, \cite{l}), which will be used throughout this
paper.\\

{\nt}A {\it{tensor norm}} $\alpha$ is a mapping which assigns to each pair $(E, F)$
of Banach spaces a norm $\alpha(\cdot; E, F)$ on the algebraic tensor product
$E \otimes F$ (shorthand: $E {\otimes}_\alpha F$  and $E \tilde{\otimes}_\alpha F$
for the completion) such that
\begin{enumerate}
\item[(1)]   $\varepsilon\leq\alpha\leq\pi$
\item[(2)]  $\alpha$ satisfies the metric mapping property: If
             $S \in {\cal L}(E, G)$ and $T \in {\cal L}(F, H)$,\\
             then $\| S\otimes T:E\otimes_\alpha F \longrightarrow G\otimes_\alpha H\|                    
             \leq \| S\|\ \| T\|$             
\end{enumerate}                          
Well-known examples are the injective tensor norm $\varepsilon$, which is the smallest one,
and the projective tensor norm $\pi$, which is the largest one. For other important
examples we refer to \cite{df}, \cite{gl}, or \cite{l}. Each tensor norm 
$\alpha$ can be extended in two natural ways. For this, denote for given 
Banach spaces $E$ and $F$
\begin{center}
FIN$(E) : = \{M \subseteq E\mid M\in $ FIN$\}$ \hspace{0.2cm}and\hspace{0.2cm}
COFIN$(E) : = \{L\subseteq E\mid E/L \in $ FIN$\}$,       
\end{center}
where FIN stands for the class of all finite-dimensional Banach spaces.
Let $z \in E \otimes F$. Then the {\it{finite hull}}\/ $\stackrel{\rightarrow}{\alpha}$
of $\alpha$ is given by
\begin{center}
\hspace{0.2cm}$\stackrel{\rightarrow}{\alpha}(z; E, F) : = \inf\{\alpha(z; M, N) \mid M \in $ FIN$(E), 
N \in $ FIN$(F), z \in M \otimes N\}$
\end{center}
and the {\it{cofinite hull}}\/ $\stackrel{\leftarrow}{\alpha}$ of $\alpha$ is
given by
\begin{center}
\hspace{0.2cm}$\stackrel{\leftarrow}{\alpha}(z; E, F) : = \sup\{\alpha(Q_K^E \otimes Q_L^F(z); E/K, F/L) \mid K \in $ COFIN$(E), 
L \in $ COFIN$(F)\}$.
\end{center}
$\alpha$ is called {\it{finitely generated}}\/ if $\alpha = \hspace{0.1cm}\stackrel{\rightarrow}{\alpha}$,
{\it{cofinitely generated}}\/ if $\alpha = \hspace{0.1cm}\stackrel{\leftarrow}{\alpha}$ (it is
always true that $\stackrel{\leftarrow}{\alpha} \hspace{0.1cm}\leq \alpha \leq \hspace{0.1cm}\stackrel{\rightarrow}{\alpha}$).
$\alpha$ is called {\it{right-accessible}} if  
$\stackrel{\leftarrow}\alpha$$(z; M, F) = \hspace{0.1cm}\stackrel{\rightarrow}\alpha$$(z; M, F)$ for all $(M, F) \in$
FIN $\times$ BAN, {\it{left-accessible}} if 
$\stackrel{\leftarrow}{\alpha}$$(z; E, N) = \hspace{0.1cm}\stackrel{\rightarrow}{\alpha}$$(z; E, N)$ for all $(E, N) \in$
BAN $\times$ FIN, and {\it{accessible}} if it is right- and left-accessible.
$\alpha$ is called {\it{totally accessible}} if $\stackrel{\leftarrow}{\alpha}
\hspace{0.1cm}= \hspace{0.1cm}\stackrel{\rightarrow}{\alpha}$.\\ 
The injective norm $\varepsilon$ is totally
accessible, the projective norm $\pi$ is accessible - but not totally accessible,
and {\sc{Pisier's}} counterexample implies the existence of a (finitely generated)
tensor norm which is neither left- nor right-accessible (see \cite{df}, 31.6).\\
There exists a powerful one-to-one correspondence between finitely generated
tensor norms and maximal Banach ideals which links thinking in terms of 
operators with "tensorial" thinking and which allows to transfer notions in
the "tensor-language" to the "operator-language" and conversely. We refer the 
reader to \cite{df} and \cite{oe1} for detailed informations
concerning this subject. Let $E, F$ be Banach spaces and $z = \sum\limits_{i = 1}^n a_i \otimes y_i$\ be
an Element in $E' \otimes F$. Then $T_z(x): = \sum\limits_{i = 1}^n < x, a_i >  y_i $
defines a finite operator $T_z \in \oid{F}{}{}(E, F)$ which is independent of the
representation of $z$ in $E' \otimes F$. 
Let $\alpha$ be a finitely generated tensor norm and (\oid{A}{}{}, \idn{A}{}{})
be a maximal Banach ideal. $\alpha$ and (\oid{A}{}{}, \idn{A}{}{}) are said to
be {\it{associated}}, notation:
\begin{center}
$(\oid{A}{}{}, \idn{A}{}{}) \sim \alpha$ \hspace{0.15cm}(shorthand: $\oid{A}{}{} \sim \alpha$, resp. 
\hspace{0.1cm}$\alpha \sim \oid{A}{}{}$)
\end{center}
if for all $M, N \in $ FIN 
\begin{center}
$\oid{A}{}{}(M, N) = M'{\otimes}_\alpha N$
\end{center}
holds isometrically: $\idn{A}{}{}(T_z) = \alpha(z; M', N)$.\\

{\nt}Important examples of Banach ideals are given by  
$(\cal I, \bf I) \sim \pi$ (integral operators),
$(\oid{L}{2}{}, \idn{L}{2}{}) \sim w_2 $  
(operators which factor through a Hilbert space),
$(\oid{D}{2}{}, \idn{D}{2}{}) \stackrel{1}{=} 
(\oid{L}{2}{\displaystyle\ast}, \idn{L}{2}{\displaystyle\ast}) \sim w_2^\ast $\/
($2$-dominated operators),
$(\oid{P}{p}{}, \idn{P}{p}{}) \sim g_p\backslash = g_q^\ast $ (absolutely $p$-summing operators),
$1 \leq p \leq \infty, \frac{1}{p} + \frac{1}{q} = 1$,
$(\oid{L}{\infty}{}, \idn{L}{\infty}{}) \stackrel{1}{=} 
(\oid{P}{1}{\displaystyle\ast}, \idn{P}{1}{\displaystyle\ast}) \sim w_\infty $
and $(\oid{L}{1}{}, \idn{L}{1}{}) \stackrel{1}{=} 
(\oid{P}{1}{\displaystyle\ast d}, \idn{P}{1}{\displaystyle\ast d}) \sim w_1 $.\\

\section{Accessible conjugate operator ideals}

{\nt}Let $\alpha$ be an arbitrary finitely generated tensor norm and $(\oid{A}{}{}, 
\idn{A}{}{}) \sim \alpha$ the associated maximal Banach ideal. Accessibility
conditions of tensor norms can be transfered to the operator ideal language in the 
following sense:\\ 

{\noindent}(\oid{A}{}{}, \idn{A}{}{}) is called {\it{right-accessible}}, if for all
$(M, F) \in$ FIN $\times$ BAN, operators $T \in {\cal L}(M, F)$ and $ \varepsilon > 0$ 
there are $N \in$ FIN$(F)$ and $S \in {\cal L}(M, N)$ such that the following
diagram commutes:\\

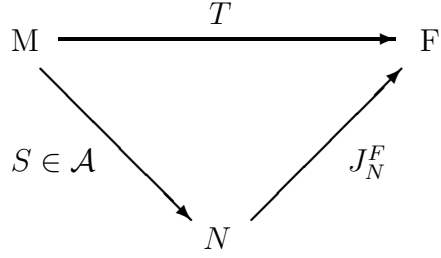
\setlength{\unitlength}{0.8cm}
\begin{picture}(6,1)\thicklines
 \put(6.8,0){\vector(1,0){5.6}}
 \put(6,-0.2){M}
 \put(12.8,-0.2){F}
 \put(6.5,-0.5){\vector(1,-1){2.5}}
 \put(9.2,-3.5){$N$}
 \put(10,-3){\vector(1,1){2.5}}
 \put(9.3,0.25){$T$}
 \put(11.6,-2.2){$J_N^F$}
 \put(6.0,-2.2){$S \in {\cal A}$}
\end{picture}\\
\\
\\
\\
\\
\\
\\
\\
and such that $\idn{A}{}{}(S) \leq (1 + \varepsilon) \cdot \idn{A}{}{}(T)$.\\

{\noindent}$(\oid{A}{}{}, \idn{A}{}{})$ is called {\it{left-accessible}}, if for all $(E, N) \in$ BAN $\times$ FIN,
operators $T \in {\cal L}(E, N)$ and $\varepsilon > 0$ there are $K \in$ COFIN$(E)$
and $S \in {\cal L}(E/K, N)$ such that the following diagram commutes

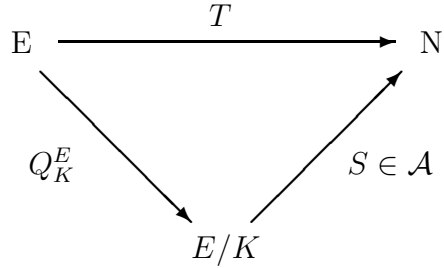
\setlength{\unitlength}{0.8cm}
\begin{picture}(6,6)\thicklines
\put(6.8,4){\vector(1,0){5.6}}     
 \put(6,3.8){E}
 \put(12.8,3.8){N}
 \put(6.5,3.5){\vector(1,-1){2.5}}
 \put(9.0,0.4){$E/K$}
 \put(10,1){\vector(1,1){2.5}}
 \put(9.3,4.25){$T$}
 \put(11.6,1.8){$S \in {\cal A}$}
 \put(6.3,1.8){$Q_K^E$}
 
\end{picture}\\

{\noindent}and such that $\idn{A}{}{}(S) \leq (1 + \varepsilon) \cdot \idn{A}{}{}(T)$.\\

{\noindent}$(\oid{A}{}{}, \idn{A}{}{})$ is {\it{totally accessible}}, if for every finite rank operator
$T \in {\cal F}(E, F)$ between arbitrary Banach spaces $E, F$ and $\varepsilon > 0$ there are
$(K, N) \in$ COFIN$(E) \times$ FIN$(F)$ and $S \in {\cal L}(E/K, N)$ such that 
the following diagram commutes

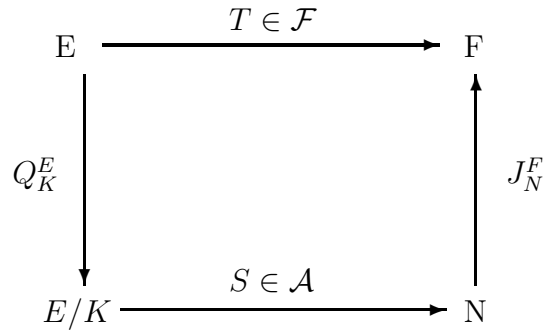
\setlength{\unitlength}{0.8cm}
\begin{picture}(6,6)\thicklines
 \put(6.8,4){\vector(1,0){5.6}}
 \put(6,3.8){E}
 \put(12.8,3.8){F}
 \put(6.5,3.5){\vector(0,-1){3.5}}
 \put(5.8,-0.6){$E/K$}
 \put(13,0){\vector(0,1){3.5}}
 \put(12.8,-0.6){N}
 \put(7.1,-0.4){\vector(1,0){5.4}}
 \put(8.9,4.25){$T \in {\cal F}$}
 \put(13.5,1.7){$J_N^F$}
\put(8.9,-0.1){$S \in {\cal A}$}
 \put(5.3,1.7){$Q_K^E$}
 
\end{picture}\\
\\
\\
\\
{\noindent}and such that $\idn{A}{}{}(S) \leq (1 + \varepsilon) \cdot \idn{A}{}{}(T)$.\\

{\noindent}The problem whether each {\it{maximal}} Banach ideal is accessible, 
was negatively answered by  {\sc{G. Pisier}} at Oberwolfach in 1991 
(cf. \cite{df}, 31.6): 
 
\begin{theorem} {\bf{(G. Pisier, 1991)}} There exists a maximal 
Banach ideal, which is not accessible.
\end{theorem}

{\noindent}Let $(\oid{A}{}{}, \idn{A}{}{})$ be a given {\it{Banach}} ideal.
Looking at the following "increasing sequence",

\begin{center}
$(\oid{A}{}{min}, \idn{A}{}{min}) \subseteq (\oid{A}{}{{\displaystyle\ast}{\scriptstyle\triangle}}, 
\idn{A}{}{{\displaystyle\ast}{\scriptstyle\triangle}}) \subseteq 
(\oid{A}{}{max}, \idn{A}{}{max})$
\end{center}

{\nt}then it follows that:

\begin{itemize} 
\item $(\oid{A}{}{min}, \idn{A}{}{min})$ is accessible and in 
general not totally accessible.
\item $(\oid{A}{}{{\displaystyle\ast}{\scriptstyle\triangle}}, 
\idn{A}{}{{\displaystyle\ast}{\scriptstyle\triangle}})$
is {\it{right}}-accessible. In particular 
$(\oid{A}{}{{\displaystyle\ast}{\scriptstyle\triangle dd}}, 
\idn{A}{}{{\displaystyle\ast}{\scriptstyle\triangle dd}})$ 
is accessible (cf. \cite{oe1}, \cite{oe2}).
\item $(\oid{A}{}{\displaystyle\ast\ast}, \idn{A}{}{\displaystyle\ast\ast})
      = (\oid{A}{}{max}, \idn{A}{}{max})$ 
      in general is not accessible.
\end{itemize}

{\noindent}{\Large}Hence, the "larger" the ideal, the "fewer" accessible is it.\\  
 
{\noindent}{\bf{OPEN PROBLEM:}} {\it{Is 
$(\oid{A}{}{{\displaystyle\ast}{\scriptstyle\triangle}}, 
\idn{A}{}{{\displaystyle\ast}{\scriptstyle\triangle}})$ always 
{\it{left}}-accessible?}}\\

{\noindent}Note, that if (\oid{A}{}{\displaystyle\ast}, \idn{A}{}{\displaystyle\ast})
is right-accessible, then 
(\oid{A}{}{\displaystyle\ast\ast}, \idn{A}{}{\displaystyle\ast\ast})
is left-accessible. In particular 
(\oid{A}{}{{\displaystyle\ast}{\scriptstyle\triangle}}, 
\idn{A}{}{{\displaystyle\ast}{\scriptstyle\triangle}})
is left-accessible. Hence, if there exists a non-left-accessible
conjugate of a maximal Banach ideal, then this 
ideal is not right-accessible. To solve this difficult problem, 
we look for conditions which are {\it{equivalent}} to the left-accessibility 
of $(\oid{A}{}{{\displaystyle\ast}{\scriptstyle\triangle}}, 
\idn{A}{}{{\displaystyle\ast}{\scriptstyle\triangle}})$. To this end 
remember the
 
\begin{theorem}{\bf{(Principle of local reflexivity)}} Let $M$ and $F$ be Banach
spaces, $M$ finite-dimensional and $T \in \oid{L}{}{}(M,F'')$. Then for 
every $\epsilon > 0$ and $N \in FIN(F')$ there exists an operator $S \in 
\oid{L}{}{}(M,F)$ such that
\begin{itemize}
\item[$(1)$]   $\|S\| \leq (1+\epsilon) \cdot \|T\|$
\item[$(2)$]   $<Sx,b> = <b,Tx>$ for all $(x,b) \in M \times N$
\item[$(3)$]   ${j_F}Sx = Tx$ for all $x \in M \cap T^{-1}(R({j_F}))$
\end{itemize}
\end{theorem}

{\nt}A transfer of the principle of local reflexivity (from the 
classical operator norm) to arbitrary $p$-norms of operator ideals,
which is directly related to the left-accessibility of conjugate 
operator ideals, is given in the following sense (cf. \cite{oe1} and \cite{oe2}):\\

{\nt}{\bf{Definition 4.1:}} Let $M$ and $F$ be Banach
spaces, $M$ finite-dimensional, $N \in FIN(F')$ and $T \in \oid{L}{}{}(M,F'')$. 
Let (\oid{A}{}{}, \idn{A}{}{}) be an arbitrary
$p$-Banach ideal ($0 < p \leq 1$) and $\epsilon > 0$. We say that the {\it{principle
of $\oid{A}{}{}$-local reflexivity}} (short: $\oid{A}{}{}-LRP$) holds, 
if there exists an operator $S \in \oid{L}{}{}(M,F)$ such that
\begin{itemize}
\item[$(1)$]    $\idn{A}{}{}(S) \leq (1+\epsilon) \cdot \idn{A}{}{\displaystyle\ast\ast}(T)$
\item[$(2)$]   $<Sx,b> = <b,Tx>$ for all $(x,b) \in M \times N$
\item[$(3)$]  ${j_F}Sx = Tx$ for all $x \in M \cap T^{-1}(R({j_F}))$
\end{itemize}

{\noindent}Using an analogous proof as in \cite{p1}, ch. 28,
it can be shown that for any $p$-Banach ideal \\(\oid{A}{}{}, \idn{A}{}{}),
($0 < p\leq 1$), the \oid{A}{}{}-LRP is already satisfied if only 
the conditions (1) and (2) of the previous definition are assumed (cf. \cite{oe2}).
In which sense does this reflect accessibility conditions? The 
answer (cf. \cite{oe2}) is given by the following 

\begin{theorem} Let $(\oid{A}{}{}, \idn{A}{}{})$ be an 
arbitrary $p$-Banach ideal $(0 < p \leq 1)$. Then the following statements
are eqivalent:

\begin{itemize}
\item[$(1)$]     $(\oid{A}{}{\scriptstyle\triangle}, 
               \idn{A}{}{\scriptstyle\triangle})$ is left-accessible

\item[$(2)$]    $\oid{A}{}{\displaystyle\ast\ast}(M, F'') \cong
               \oid{A}{}{}(M, F)^{''}$ for all $(M, F) \in FIN 
               \times BAN$

\item[$(3)$]   The \oid{A}{}{}-LRP holds. 
\end{itemize}

\end{theorem}

{\nt}To obtain operator ideals which satisfy the transfer of the norm
estimation in the \oid{L}{}{}-LRP to their ideal norm, we need further 
geometrical properties of such operators (for an interesting connection 
of the \oid{A}{}{}-LRP for injective Banach ideals \oid{A}{}{} with
Grothendieck's inequality we refer the reader to \cite{oe4}). First let us 
note the following

\begin{theorem} Let $(\oid{A}{}{}, \idn{A}{}{})$ be an arbitrary 
Banach ideal. If the \oid{A}{}{\displaystyle\ast}-LRP holds, then 
$\oid{A}{}{inj} \circ \bar{\oid{F}{}{}}$ is totally accessible.
\end{theorem}

{\nt}{\sc{PROOF}}:{\hspace{0.25cm}}Let $E, F$ be arbitrary Banach spaces and 
$L \in \oid{F}{}{}(E,F)$ an arbitrary finite operator. We put $\oid{B}{}{} : = 
\oid{A}{}{{\displaystyle\ast}{\scriptstyle\triangle}}$. Without further 
assumptions on the Banach spaces (such as approximation 
properties), the inclusion $({\cal A}^{inj})^{min}(E, F) 
\stackrel{1}{\subseteq} ({\cal A}^{min})^{inj}(E, F)$ is always true 
(cf. \cite{df}, 25.11). Since $\oid{A}{}{min} \stackrel{1}{\subseteq} \oid{B}{}{}$, 
it follows therefore 
that $({\cal A}^{inj})^{min} \stackrel{1}{\subseteq} ({\cal A}^{min})^{inj}
\stackrel{1}{\subseteq} \oid{B}{}{inj}$. By theorem 4.3 and the asssumption,
\oid{B}{}{} is left-accessible, which implies that \oid{B}{}{inj} 
is totally accessible. Hence, given $\epsilon > 0$, there exist Banach 
spaces $K \in$ COFIN($E$), $N \in$ FIN($F$) and an operator $A \in 
\oid{L}{}{}(E/K, N)$, such that $L = J_{N}^{F} A Q_{K}^{E}$ and 
\begin{center}
$({\bf A}^{inj})^{min}(L) 
\stackrel{1}{=} ({\bf B}^{inj})^{min}(L) \leq \idn{B}{}{inj}(A) 
< (1 + \epsilon) \cdot \idn{B}{}{inj}(L)  
\leq (1 + \epsilon) \cdot ({\bf A}^{inj})^{min}(L)$. 
\end{center}
Therefore, for {\it{all finite operators}} we have obtained the following 
identities:
\begin{center}
$({\bf A}^{inj})^{min}(L) = \idn{B}{}{inj}(L) = ({\bf A}^{min})^{inj}(L).$
\end{center}   
Since $({\cal A}^{inj})^{min} \stackrel{1}{=} \oid{A}{}{inj} \circ 
\bar{\oid{F}{}{}}$ (cf. \cite{df}, 25.2), the proof is finished.
${}_{\displaystyle\Box}$\\

{\nt}This theorem leads to the conjecture that \oid{A}{}{inj} even is totally 
accessible, if we "only" assume that the $\oid{A}{}{\displaystyle\ast}$-LRP
holds; this conjecture remains still open.  However, an additional extension 
property leads to further interesting aspects concerning relations between 
accessibility conditions and the local reflexivity principle for operator 
ideals.\\

{\nt}{\bf{Definition 4.2:}} Let (\oid{A}{}{}, \idn{A}{}{}) be a $p$-Banach
ideal $(0 < p \leq 1)$. We say that the $\oid{A}{}{}$-{\it{extension property}} (short: 
$\oid{A}{}{}$-{\it{EP}}) holds, if for every $\epsilon > 0$, for every metric
injection $J : E \stackrel{1}{\hookrightarrow} G$ and $T \in \oid{A}{}{}(E, F)$ 
there exists a $\tilde{T} \in \oid{A}{}{}(G, F)$ such that 
$T = \tilde{T}J$ and $\idn{A}{}{}(\tilde{T}) \leq (1+\epsilon) \cdot 
\idn{A}{}{}(T)$.\\

{\nt}One example of a maximal and injective Banach ideal for which this 
extension property holds, is given by the class of all absolutely 2-summing 
operators: $(\oid{P}{2}{}, \idn{P}{2}{}) = (\oid{P}{2}{\displaystyle\ast}, 
\idn{P}{2}{\displaystyle\ast})$ satisfies the $\oid{P}{2}{} \stackrel{1}{=} 
\oid{P}{2}{\displaystyle\ast}$-EP. 
This follows immediately by a well known factorization theorem for absolutely
2-summing operators cf. (\cite{p1}, 17.3.7) and the metric extension property 
of spaces of type $C(K)$, where K is a compact space.\\
{\nt}Further examples are given by certain minimal Banach ideals: Let 
$\oid{A}{}{} \sim \alpha$ be associated. Then 
$\oid{A}{}{inj \displaystyle\ast} \stackrel{1}{=} 
\setminus\oid{A}{}{\displaystyle\ast} \sim \setminus\alpha^{\displaystyle\ast}$.
By the representation theorem for minimal operator ideals the canonical map 
$E' \tilde{\otimes}_{\setminus\alpha^{\displaystyle\ast}} F
\rightarrow (\setminus\oid{A}{}{\displaystyle\ast})^{min}(E, F)$ 
is a metric surjection for {\it{all}} Banach spaces $E$ and $F$ 
(cf. \cite{df}, 22.2). Since $\setminus\alpha^{\displaystyle\ast}$ is 
left-projective, it follows therefore that the 
$(\setminus\oid{A}{}{\displaystyle\ast})^{min}$-EP always holds.\\
Note that the \oid{L}{}{}-EP is false, since there is no Hahn-Banach theorem
for operators.\\

{\nt}To prepare the proof of the following theorem, observe that for a given 
$p$-Banach ideal $(\oid{A}{}{}, \idn{A}{}{})$ we have the inclusion 
$\oid{A}{}{} \stackrel{1}{\subseteq} 
\oid{A}{}{\scriptstyle\triangle \scriptstyle\triangle}$ on the class of all 
{\it{finite}} operators and the (global) inclusion 
$\oid{A}{}{\scriptstyle\triangle} \stackrel{1}{\subseteq} 
\oid{A}{}{\displaystyle\ast}$ which implies that 
$\oid{A}{}{\displaystyle\ast \scriptstyle\triangle} 
\stackrel{1}{\subseteq} \oid{A}{}{\scriptstyle\triangle \scriptstyle\triangle}$.  
These inclusions follow directly by the definition of conjugate and adjoint 
operator ideals. Since there exists a maximal Banach ideal which is not 
right-accessible and since \oid{A}{}{\displaystyle\ast \scriptstyle\triangle}
always is right-accessible, it follows that in general  
$\oid{A}{}{\displaystyle\ast \scriptstyle\triangle} 
\stackrel{1}{\not=} \oid{A}{}{\scriptstyle\triangle \scriptstyle\triangle}$.

\begin{theorem} Let $(\oid{A}{}{}, \idn{A}{}{})$ be  an arbitrary $p$-Banach
ideal $(0 < p \leq 1)$.
\begin{itemize}
\item[(1)] If the \oid{A}{}{}-EP holds,  then \oid{A}{}{} is left-accessible
           and $\oid{A}{}{} \stackrel{1}{\subseteq} 
           \oid{A}{}{\scriptstyle\triangle \scriptstyle\triangle}$. 
           In particular \oid{A}{}{inj} is totally accessible.  
\item[(2)] If the \oid{A}{}{\displaystyle\ast}-EP and the 
           \oid{A}{}{\scriptstyle\triangle}-LRP both are given, then  
           \oid{A}{}{} is accessible and 
           $\oid{A}{}{\scriptstyle\triangle \scriptstyle\triangle} 
           \stackrel{1}{=} \oid{A}{}{\displaystyle\ast \scriptstyle\triangle}$.  
\end{itemize} 
\end{theorem}

{\nt}{\sc{PROOF}}:{\hspace{0.25cm}}To prove (1), let $T \in \oid{F}{}{}(E, N)$
an arbitrary (finite) operator, considered as an element of $\oid{A}{}{}(E, N)$,
where $(E, N) \in $ BAN $\times$ FIN. Let 
$J_E : E \stackrel{1}{\hookrightarrow} E^{\infty}$ the canonical injection 
from the Banach space $E$ into the Banach space $E^{\infty} = 
C(B_{E'})$. Given $\epsilon > 0$, the assumption of the \oid{A}{}{}-EP
implies the existence of an operator $\tilde{T} \in \oid{A}{}{}(E^{\infty}, N)$, 
such that $T = \tilde{T} J_{E}$ and $\idn{A}{}{}(\tilde{T}) \leq (1 + \epsilon)
\cdot \idn{A}{}{}(T)$. Since the dual of $E^{\infty}$ has the metric approximation 
property, it follows that there exists a {\it{finite}} operator $A \in 
\oid{F}{}{}(E^{\infty}, E^{\infty})$, such that $\tilde{T} = \tilde{T} A$ and
$\| A \| < 1 + \epsilon$. Hence we obtain that $T = Id_N \tilde{T} A J_E
\in \oid{A}{}{min}(E, N)$ and
\begin{center}
${\bf A}^{min}(T) = {\bf A}^{min}(Id_N \tilde{T} A J_E) \leq 
\idn{A}{}{}(\tilde{T}) \cdot \| A J_{E} \| \leq (1 + \epsilon)^{2} \cdot 
\idn{A}{}{}(T) \leq (1 + \epsilon)^{2} \cdot \idn{A}{}{min}(T)$. 
\end{center}
Since minimal $p$-Banach ideals are always accessible, the last estimation 
shows that \oid{A}{}{} is left-accessible. A straightforward calculation  
shows that $\oid{A}{}{}(E^{\infty}, F) \stackrel{1}{\subseteq} 
\oid{A}{}{\scriptstyle\triangle \scriptstyle\triangle}(E^{\infty}, F)$ for all 
Banach spaces $E$ and $F$. Hence, the \oid{A}{}{}-EP implies that 
$\oid{A}{}{} \stackrel{1}{\subseteq} 
\oid{A}{}{\scriptstyle\triangle \scriptstyle\triangle}$. \\ 

{\nt}To prove statement (2), let $E, F$ be arbitrary Banach spaces, $\epsilon 
> 0$ and set 
$\oid{B}{}{} : = \oid{A}{}{{\displaystyle\ast}{\scriptstyle\triangle}}$.
First we show that 
\begin{itemize} 
\item[$(*)$]  \hspace{4.0 cm}$\idn{B}{}{}(T'') \leq \idn{A}{}{}(T)$ for all 
            $T \in \oid{F}{}{}(E, F)$. 
\end{itemize}
To this end let $T \in \oid{F}{}{}(E, F)$ and $L \in \oid{F}{}{}(F'', E'')$ 
be arbitrary finite operators. By the assumed $\oid{A}{}{\displaystyle\ast}$-EP
there exists an operator $\tilde{L} \in 
\oid{A}{}{\displaystyle\ast}((F'')^{\infty}, E'')$ such that $L = \tilde{L} J_{F''}$
and $\idn{A}{}{\displaystyle\ast}(\tilde{L}) \leq (1 + \epsilon) \cdot 
\idn{A}{}{\displaystyle\ast}(L)$. As in the proof of (1), we find an      
operator $A \in \oid{F}{}{}((F'')^{\infty}, (F'')^{\infty})$, 
such that $T'' \tilde{L} = T'' \tilde{L} A$ and $\| A \| \leq 1 + \epsilon$. Canonical
factorization of $A$ leads to a finite dimensional Banach space $N$, operators 
$U \in \oid{L}{}{}(N, (F'')^{\infty})$ and $V \in 
\oid{L}{}{}((F'')^{\infty}, N)$ with $A = U V$, $\| V \| \leq  1$ and 
$\| U \|  \leq  \| A \| \leq 1 + \epsilon $. Since $\tilde{L} U \in 
\oid{F}{}{}(N, E'')$, the assumed validity of the 
\oid{A}{}{\scriptstyle\triangle}-LRP implies the existence of an operator 
$S \in \oid{L}{}{}(N, E)$ such that $T'' \tilde{L} U = T'' j_E S = 
j_F T S$ and $\idn{A}{}{\scriptstyle\triangle}(S) \leq (1 + \epsilon) 
\cdot \idn{A}{}{\displaystyle\ast}(\tilde{L} U)$ \/ (cf. \cite{oe2}, lemma 1.1).
Hence $T'' L = T'' \tilde{L} J_{F''} = j_F T S V J_{F''}$, which implies
that:
\begin{center}
$\mid tr(T'' L) \mid = \mid tr(j_F T S V J_{F''}) \mid = 
\mid tr(S V J_{F''} j_F T) \mid$ 
\end{center}
Since $T$ is a finite operator, we therefore obtain the following estimation:
\begin{center}
$\mid tr(T'' L) \mid \leq \idn{A}{}{\scriptstyle\triangle}(S V) \cdot \idn{A}{}{}(T) 
\leq (1 + \epsilon) \cdot \idn{A}{}{\displaystyle\ast}(\tilde{L} U) 
\cdot \idn{A}{}{}(T) \leq (1 + \epsilon)^3 \cdot \idn{A}{}{\displaystyle\ast}(L) 
\cdot \idn{A}{}{}(T),$
\end{center}
which implies $(*)$. Obviously $\oid{B}{}{dd}
\stackrel{1}{\subseteq} \oid{B}{}{}$, and therfore it follows that
$\oid{A}{}{} \stackrel{1}{\subseteq} \oid{B}{}{dd} \stackrel{1}{\subseteq}
\oid{B}{}{}$ on the class of all {\it{finite}} operators. Conjugation
of this inclusion implies $\oid{B}{}{\scriptstyle\triangle} \stackrel{1}{\subseteq} 
\oid{A}{}{\scriptstyle\triangle}$, and a further conjugation leads to
$\oid{A}{}{\scriptstyle\triangle \scriptstyle\triangle} \stackrel{1}{\subseteq}
\oid{B}{}{\scriptstyle\triangle \scriptstyle\triangle} \stackrel{1}{=}
\oid{B}{}{} \stackrel{1}{\subseteq}
\oid{A}{}{\scriptstyle\triangle \scriptstyle\triangle}$. Hence
$\oid{A}{}{\scriptstyle\triangle \scriptstyle\triangle} \stackrel{1}{=}
\oid{A}{}{\displaystyle\ast \scriptstyle\triangle}$. 
Obviously, the \oid{A}{}{\scriptstyle\triangle}-LRP implies the 
\oid{A}{}{\displaystyle\ast}-LRP. Hence, $\oid{B}{}{} \stackrel{1}{=}
\oid{A}{}{\scriptstyle\triangle \scriptstyle\triangle}$
is accessible, which immediately leads to the accessibility of \oid{A}{}{}.
${}_{\displaystyle\Box}$

\begin{cor}
Let $(\oid{A}{}{}, \idn{A}{}{})$ be a maximal Banach ideal with the 
\oid{A}{}{}-EP. Then the \oid{A}{}{\displaystyle\ast}-LRP
holds.
\end{cor}

{\nt}Recall the following (cf. \cite{oe1}, Satz 3.10):

\begin{theorem} Let $(\oid{A}{}{}, \idn{A}{}{})$ 
be an arbitrary left-accessible Banach ideal. Then the following statements 
are equivalent:
\begin{itemize}
\item[(i)]     the \oid{A}{}{\scriptstyle\triangle}-LRP holds 
\item[(ii)]    the \oid{A}{}{\displaystyle\ast}-LRP holds .
\end{itemize}
\end{theorem}

{\nt}Combining the last two theorems, another interesting fact follows:

\begin{theorem} Let $(\oid{A}{}{}, \idn{A}{}{})$ be an arbitrary 
maximal and left-accessible Banach ideal, such that the 
\oid{A}{}{\displaystyle\ast}-EP holds. Then \oid{A}{}{} is accessible and 
$\oid{A}{}{\scriptstyle\triangle \scriptstyle\triangle} \stackrel{1}{=} 
\oid{A}{}{\displaystyle\ast \scriptstyle\triangle}$.  
\end{theorem}

{\nt}The previous considerations naturally lead to the following (still unsolved) 
problem: Is it possible to drop the assumption "left-accessible" in the 
theorems 4.6 and 4.7 ?

\
\\
\\
\\
\\
{\nt}{\it{Concerning further discussions, please contact:}}\\ 

{\noindent}Frank Oertel\\
Swiss Reinsurance Company\\
Dpt.: PM-PH (Development and Support)\\
Mythenquai 50/60\\
CH - 8022 Zurich\\
SWITZERLAND\\

{\noindent}E-mail: frank.oertel@swissre.ch\\
Tel.: +41-1-285-3688\\

\end{document}